\documentclass[11pt,a4paper,reqno]{amsart}  
\usepackage{amssymb} \usepackage{amscd} \usepackage{times}
\newcommand{\be} {\begin{equation}}
\newcommand{\ee} {\end{equation}}
\newcommand{\bea} {\begin{eqnarray}}
\newcommand{\eea} {\end{eqnarray}}
\newcommand{\Bea} {\begin{eqnarray*}}
\newcommand{\Eea} {\end{eqnarray*}}

\def\zbb{\mathbb{Z}}  
  
  \def\phi{\varphi}
 \def\p1{{\mathbb{P}^1_\zbb}}

\newtheorem{Theorem}{\quad Theorem}[section]

\begin{document}
\title{ A compactness result for an elliptic equation in dimension 2.}
\author{Samy Skander Bahoura} 
\address{Departement de Mathematiques, Universite Pierre et Marie Curie, 2 place Jussieu, 75005, Paris, France.}
\email{samybahoura@yahoo.fr, samybahoura@gmail.com}
\maketitle

\begin{abstract}
We give a blow-up behavior for the solutions of an elliptic equation under some conditions. We also derive a compactness creterion for this equation.
\end{abstract}



\section{Introduction and Main Results}

Let us consider the following operator:
$$ L_{\epsilon}:=\Delta +\epsilon( x_1 \partial_1+ x_2 \partial_2) =\dfrac{div[a_{\epsilon}(x) \nabla]}{a_{\epsilon}(x)}\,\, {\rm with} \,\, a_{\epsilon}(x)= e^{\frac{\epsilon |x|^2}{2}}. $$
We consider the following equation:
\begin{displaymath}  (P_{\epsilon})  \left \{ \begin {split} 
      -\Delta u-\epsilon( x_1 \partial_1 u+ x_2 \partial_2 u)=-L_{\epsilon} u & = V e^{u}     \,\, &&\text{in} \!\!&&\Omega \subset {\mathbb R}^2, \\
                                                                                                                u  & = 0  \,\,            && \text{in} \!\!&&\partial \Omega.               
\end {split}\right.
\end{displaymath}
Here, we assume that:
$$ \Omega \,\, \text{starshaped}, $$
and,
$$ u \in W_0^{1,1}(\Omega), \,\,  e^u \in L^1(\Omega),  \,\, 0 \leq V \leq b, \,\, 1\geq \epsilon \geq 0 . $$

When $ \epsilon =0 $ the previous equation was studied by many authors with or without  the boundary condition, also for Riemann surfaces see [1-20] where one can find some existence and compactness results. Also we have a nice formulation in the sens of the distributions of this Problem in [7].

Among other results, we  can see in [6] the following important Theorem,

{\bf Theorem A }{\it (Brezis-Merle [6])} {\it If $ (u_i)_i $ and $ (V_i)_i $ are two sequences of functions relative to the problem $ (P_0) $ with $ \epsilon =0 $ and,
$$ 0 < a \leq V_i \leq b < + \infty $$
then it holds,
$$ \sup_K u_i \leq c, $$
with $ c $ depending on $ a, b,  K $ and $ \Omega $.}

We can find in [6] an interior estimate if we assume $ a=0 $ but we need an assumption on the integral of $ e^{u_i} $, namely:

{\bf Theorem B}{\it (Brezis-Merle [6])}.{\it For $ (u_i)_i $ and $ (V_i)_i $ two sequences of functions relative to the problem $ (P_0) $ with,
$$ 0 \leq V_i \leq b < + \infty \,\, {\rm and} \,\, \int_{\Omega} e^{u_i} dy  \leq C, $$
then it holds;
$$ \sup_K u_i \leq c, $$
with $ c $ depending on $ b, C, K $ and $ \Omega $.}

The condition $ \int_{\Omega} e^{u_i} dy  \leq C $ is a necessary condition in the Problem $ (P_{\epsilon}) $ as showed by the following counterexample for $ \epsilon =0 $:

{\bf Theorem C} {\it (Brezis-Merle [6])}.{\it There are two sequences $ (u_i)_i $ and $ (V_i)_i $ of the problem $ (P_0) $ with;
$$ 0 \leq V_i \leq b < + \infty, \, \, \,  \int_{\Omega} e^{u_i} dy  \leq C, $$
such that,
$$ \sup_{\Omega}  u_i \to + \infty. $$}

To obtain the two first previous results (Theorems A and B) Brezis and Merle used  an inequality (Theorem 1 of [6]) obtained by an approximation argument and they used Fatou's lemma and applied the maximum principle in $ W_0^{1,1}(\Omega) $ which arises from Kato's inequality. Also this weak form of the maximum principle is used to prove the local uniform boundedness result by comparing  a certain function and the Newtonian potential. We refer to [5] for a topic about the weak form of the maximum principle.

Note that for the problem $ (P_0) $, by using the Pohozaev identity, we can prove that $ \int_{\Omega} e^{u_i} $ is uniformly bounded when $ 0 < a \leq V_i \leq b < +\infty $ and $  ||\nabla V_i||_{L^{\infty}} \leq A $ and $ \Omega $ starshaped, when $ a=0 $ and $ \nabla \log V_i $ is uniformly bounded, we can bound uniformly $ \int_{\Omega} V_i e^{u_i} $  . In [17] Ma-Wei have proved that those results stay true for all open sets not necessarily starshaped in the case $ a>0 $.

In [8] Chen-Li have proved that if $ a=0 $ and $ \int_{\Omega} e^{u_i} $ is uniformly bounded and $ \nabla \log V_i $ is uniformly bounded then $ (u_i)_i $ is bounded near the boundary and  we have directly the compactness result for the problem $ (P_0) $. Ma-Wei in [17] extend this result in the case where $ a >0 $.

When $ \epsilon =0 $  and if we assume $ V $ more regular we can have another type of estimates called $ \sup + \inf $ type inequalities. It was proved by Shafrir see [19] that, if $ (u_i)_i, (V_i)_i $ are two sequences of functions solutions of  the Problem $ (P_0) $ without assumption on the boundary and $ 0 < a \leq V_i \leq b < + \infty $  then it holds:

$$ C\left (\dfrac{a}{b} \right ) \sup_K u_i + \inf_{\Omega} u_i \leq c=c(a, b, K, \Omega). $$

We can see in [9]  an explicit value of $ C\left (\dfrac{a}{b}\right ) =\sqrt {\dfrac{a}{b}} $. In his proof, Shafrir has used the blow-up function, the Stokes formula and an isoperimetric inequality see [2]. For Chen-Lin, they have used the blow-up analysis combined with some geometric type inequality for the integral curvature see [9].

Now, if we suppose $ (V_i)_i $ uniformly Lipschitzian with $ A $ its Lipschitz constant then $ C(a/b)=1 $ and $ c=c(a, b, A, K, \Omega)
$ see Brezis-Li-Shafrir [4]. This result was extended for
H\"olderian sequences $ (V_i)_i $ by Chen-Lin see  [9]. Also have in [15], an extension of the Brezis-Li-Shafrir result to compact
Riemannian surfaces without boundary. One can see in [16] explicit form,
($ 8 \pi m, m\in {\mathbb N}^* $ exactly), for the numbers in front of
the Dirac masses when the solutions blow-up. Here the notion of isolated blow-up point is used. Also one can see in [10] refined estimates near the isolated blow-up points and the  bubbling behavior  of the blow-up sequences.

Here we give the behavior of the blow-up points on the boundary and a proof of a compactness result with Lipschitz condition. Note that our problem is an extension of the Brezis-Merle Problem.

The Brezis-Merle Problem (see [6]) is:

{\bf Problem.} Suppose that $ V_i \to  V $ in $ C^0( \bar \Omega ) $ with $ 0 \leq V_i $. Also, we consider a sequence of solutions $ (u_i) $ of $ (P_0) $ relative to $ (V_i) $ such that,
$$ \int_{\Omega} e^{u_i} dx \leq C,  $$
is it possible to have:
$$ ||u_i||_{L^{\infty}}\leq C=C(b, C, V, \Omega) ? $$

Here we give blow-up analysis on the boundary  when $ V $ (similar to the prescribed curvature when $ \epsilon =0 $) are nonegative and bounded,  and on the other hand, if we add the assumption that these functions (similar to the prescribed cruvature) are uniformly Lipschitzian, we have a compactness of the solutions of the problem $ (P_{\epsilon}) $ for $ \epsilon $ small enough. (In particular we can take a sequence of $ \epsilon_i $ tending to $ 0 $):

For the behavior of the blow-up points on the boundary, the following condition is sufficient,
$$ 0 \leq  V_i \leq b, $$
The condition $ V_i \to  V $ in $ C^0(\bar \Omega) $ is not necessary. But for the compactness of the solutions we add the following condition:
$$ ||\nabla V_i||_{L^{\infty}}\leq  A. $$

Our main results are:

\begin{Theorem}  Assume that $ \max_{\Omega} u_i \to +\infty $, where $ (u_i) $ are solutions of the probleme $ (P_{\epsilon_i}) $ with:
 $$ 0 \leq V_i \leq b,\,\,\, {\rm and } \,\,\, \int_{\Omega}  e^{u_i} dx \leq C, \,\,\, \epsilon_i \to 0, $$
 then,  after passing to a subsequence, there is a finction $ u $,  there is a number $ N \in {\mathbb N} $ and $ N  $ points $ x_1, \ldots, x_N \in  \partial \Omega $, such that, 
$$ \partial_{\nu} u_i  \to \partial_{\nu} u +\sum_{j=1}^N \alpha_j \delta_{x_j}, \,\,\, \alpha_j \geq 4\pi, \,\, {\rm \,\, in \, the \, sens \, of \, measures \, on \,\, \partial \Omega.} $$
$$ u_i \to u \,\,\, {\rm in }\,\,\, C^1_{loc}(\bar \Omega-\{x_1,\ldots, x_N \}). $$
\end{Theorem} 

\begin{Theorem}Assume that $ (u_i) $ are solutions of $ (P_{\epsilon_i}) $ relative to $ (V_i) $ with the following conditions:
$$ 0 \leq V_i \leq b,\,\, ||\nabla V_i||_{L^{\infty}} \leq A \,\,\, {\rm and } \,\,\, \int_{\Omega} e^{u_i} \leq C, \,\, \epsilon_i \to 0. $$
Then we have:
$$  || u_i||_{L^{\infty}} \leq c(b, A, C, \Omega), $$
\end{Theorem} 

\section{Proof of the theorems} 

\underbar {\it Proof of theorem 1.1:} 

First remark that:

\begin{displaymath}  \left \{ \begin {split} 
      -\Delta u_i &=\epsilon_i( x_1 \partial_1 u_i+ x_2 \partial_2 u_i)+V_i e^{u_i} \,\in L^1(\Omega)    \,\, &&\text{in} \!\!&&\Omega \subset {\mathbb R}^2, \\
                                                                                                                u_i  & = 0  \,\,            && \text{in} \!\!&&\partial \Omega.               
\end {split}\right.
\end{displaymath}

and,

$$ u_i \in W_0^{1,1}(\Omega). $$

By the corollary 1 of Brezis-Merle see [6] we have $ e^{u_i} \in L^k(\Omega) $ for all $ k  >2 $ and the elliptic estimates of Agmon and the Sobolev embedding see [1] imply that:

$$ u_i \in W^{2, k}(\Omega)\cap C^{1, \epsilon}(\bar \Omega). $$

Also remark that, we have for two positive constants $ C_q= C(q,\Omega) $ and $ C_1=C_1(\Omega) $ (see [7]) :

$$ ||\nabla u_i||_{L^q} \leq C_q ||\Delta u_i||_{L^1} \leq (C'_q +\epsilon  C_1||\nabla u_i||_{L^1}), \,\,\forall \,\, i\,\, {\rm and  }  \,\, 1< q < 2. $$

Thus, if $ \epsilon >0 $ is small enough and by the Holder inequality, we have the following estimate:

$$ ||\nabla u_i||_{L^q} \leq C''_q, \,\,\forall \,\, i\,\, {\rm and  }  \,\, 1< q < 2. $$

\underbar {\it Step 1: interior estimate} 

First remark that, if we consider the following equation:

\begin{displaymath}   \left \{ \begin {split} 
      -\Delta w_i & = \epsilon_i (x_1\partial_1 u_i+ x_2\partial_2 u_i)\in L^q,\, 1<q <2     \,\, &&\text{in} \!\!&&\Omega \subset {\mathbb R}^2, \\
                w_i  & = 0  \,\,                                                                                                         && \text{in} \!\!&&\partial \Omega.               
\end {split}\right.
\end{displaymath}

If we consider  $ v _i $ the Newtonnian potential of $ \epsilon_i (x_1\partial_1 u_i+ x_2\partial_2 u_i) $, we have:

$$ v_i \in C^0(\bar \Omega),  \,\,\, \Delta (w_i-v_i)=0. $$

By the maximum principle $ w_i-v_i \in C^0(\bar \Omega) $ and thus $ w_i \in C^0(\bar \Omega) $.

Also we have by the elliptic estimates that $ w_i \in W^{2, 1+\epsilon} \subset L^{\infty} $, and we can write the equation of the Problem as:

\begin{displaymath}  \left \{ \begin {split} 
      -\Delta (u_i-w_i) &= \tilde V_i e^{u_i-w_i}   \,\, &&\text{in} \!\!&&\Omega \subset {\mathbb R}^2, \\
                   u_i-w_i  & = 0  \,\,                               && \text{in} \!\!&&\partial \Omega.               
\end {split}\right.
\end{displaymath}

with,

$$ 0 \leq \tilde V_i =V_i e^{w_i} \leq \tilde b,\, \int_{\Omega} e^{u_i-w_i} \leq \tilde C. $$

We apply the Brezis-Merle theorem to $ u_i-w_i $ to have:

$$ u_i-w_i \in L_{loc}^{\infty}(\Omega), $$

and, thus:

$$ u_i \in L_{loc}^{\infty}(\Omega). $$

\underbar {\it Step2: boundary estimate}

Set $ \partial_{\nu} u_i $ the inner derivative of $ u_i $. By the maximum principle $ \partial_{\nu} u_i \geq 0 $. 

We have:

$$ \int_{\partial \Omega} \partial_{\nu} u_i d\sigma \leq C. $$

 We have the existence of a nonnegative Radon measure $ \mu $ such that,

$$ \int_{\partial \Omega} \partial_{\nu} u_i \phi  d\sigma \to \mu(\phi), \,\,\, \forall \,\,\, \phi \in C^0(\partial \Omega). $$

We take an $ x_0 \in \partial \Omega $ such that, $ \mu({x_0}) < 4\pi $. Set $ B(x_0, \epsilon) \cap \partial \Omega := I_{\epsilon} $. We choose a function $ \eta_{\epsilon} $ such that,

$$ \begin{cases}
    
\eta_{\epsilon} \equiv 1,\,\,\,  {\rm on } \,\,\,  I_{\epsilon}, \,\,\, 0 < \epsilon < \delta/2,\\

\eta_{\epsilon} \equiv 0,\,\,\, {\rm outside} \,\,\, I_{2\epsilon }, \\

0 \leq \eta_{\epsilon} \leq 1, \\

||\nabla \eta_{\epsilon}||_{L^{\infty}(I_{2\epsilon})} \leq \dfrac{C_0(\Omega, x_0)}{\epsilon}.

\end{cases} $$

We take a $\tilde \eta_{\epsilon} $ such that,

\begin{displaymath}   \left \{ \begin {split} 
      -\Delta \tilde \eta_{\epsilon} & = 0  \,\, &&\text{in} \!\!&&\Omega \subset {\mathbb R}^2, \\
                   \tilde\eta_{\epsilon} & =  \eta_{\epsilon} \,\,                                 && \text{in} \!\!&&\partial \Omega.               
\end {split}\right.
\end{displaymath}

{\bf Remark:} We use the following steps in the construction of $ \tilde \eta_{\epsilon} $:

We take a cutoff function $ \eta_{0} $ in $ B(0, 2) $ or $ B(x_0, 2) $:

1- We set $ \eta_{\epsilon}(x)= \eta_0(|x-x_0|/\epsilon) $ in the case of the unit disk it is sufficient.

2- Or, in the general case: we use a chart $ (f, \tilde \Omega) $ with $ f(0)=x_0 $ and we take $ \mu_{\epsilon}(x)= \eta_0 ( f( |x|/ \epsilon)) $ to have  connected  sets $ I_{\epsilon} $ and we take $ \eta_{\epsilon}(y)= \mu_{\epsilon}(f^{-1}(y))$. Because $ f, f^{-1} $ are Lipschitz, $ |f(x)-x_0| \leq k_ 2|x|\leq 1 $ for $ |x| \leq 1/k_2 $ and $ |f(x)-x_0| \geq k_ 1|x|\geq 2 $ for $ |x| \geq 2/k_1>1/k_2 $, the support  of $ \eta $ is in $ I_{(2/k_1)\epsilon} $.

$$ \begin{cases}
    
\eta_{\epsilon} \equiv 1,\,\,\,  {\rm on } \,\,\,  f(I_{(1/k_2)\epsilon}), \,\,\, 0 < \epsilon < \delta/2,\\

\eta_{\epsilon} \equiv 0,\,\,\, {\rm outside} \,\,\, f(I_{(2/k_1)\epsilon }), \\

0 \leq \eta_{\epsilon} \leq 1, \\

||\nabla \eta_{\epsilon}||_{L^{\infty}(I_{(2/k_1)\epsilon})} \leq \dfrac{C_0(\Omega, x_0)}{\epsilon}.

\end{cases} $$

3- Also, we can take: $ \mu_{\epsilon}(x)= \eta_0(|x|/\epsilon) $ and $ \eta_{\epsilon}(y)= \mu_{\epsilon}(f^{-1}(y)) $, we extend it by $ 0 $ outside $ f(B_1(0)) $.  We have $ f(B_1(0)) = D_1(x_0) $, $ f (B_{\epsilon}(0))= D_{\epsilon}(x_0) $ and $ f(B_{\epsilon}^+)= D_{\epsilon}^+(x_0) $ with $ f $ and $ f^{-1} $ smooth diffeomorphism.

$$ \begin{cases}
    
\eta_{\epsilon} \equiv 1,\,\,\,  {\rm on \, a \, the \, connected \, set } \,\,\,  J_{\epsilon} =f(I_{\epsilon}), \,\,\, 0 < \epsilon < \delta/2,\\

\eta_{\epsilon} \equiv 0,\,\,\, {\rm outside} \,\,\, J'_{\epsilon}=f(I_{2\epsilon }), \\

0 \leq \eta_{\epsilon} \leq 1, \\

||\nabla \eta_{\epsilon}||_{L^{\infty}(J'_{\epsilon})} \leq \dfrac{C_0(\Omega, x_0)}{\epsilon}.

\end{cases} $$

And, $ H_1(J'_{\epsilon}) \leq C_1 H_1(I_{2\epsilon}) = C_1 4\epsilon $, because $ f $ is Lipschitz. Here $ H_1 $ is the Hausdorff measure.

We solve the Dirichlet Problem:

\begin{displaymath}  \left \{ \begin {split} 
      \Delta \bar \eta_{\epsilon}  & = \Delta \eta_{\epsilon}              \,\, &&\text{in} \!\!&&\Omega \subset {\mathbb R}^2, \\
                  \bar \eta_{\epsilon} & = 0   \,\,             && \text{in} \!\!&&\partial \Omega.               
\end {split}\right.
\end{displaymath}

and finaly we set $ \tilde \eta_{\epsilon} =-\bar \eta_{\epsilon} + \eta_{\epsilon} $. Also, by the maximum principle and the elliptic estimates we have :

$$ ||\nabla \tilde \eta_{\epsilon}||_{L^{\infty}} \leq C(|| \eta_{\epsilon}||_{L^{\infty}} +||\nabla \eta_{\epsilon}||_{L^{\infty}} + ||\Delta \eta_{\epsilon}||_{L^{\infty}}) \leq \dfrac{C_1}{\epsilon^2}, $$

with $ C_1 $ depends on $ \Omega $.

As we said in the beguening,  see also [3, 7, 13, 20], we have:

$$ ||\nabla u_i||_{L^q} \leq C_q, \,\,\forall \,\, i\,\, {\rm and  }  \,\, 1< q < 2. $$

We deduce from the last estimate that, $ (u_i) $ converge weakly in $ W_0^{1, q}(\Omega) $, almost everywhere to a function $ u \geq 0 $ and $ \int_{\Omega} e^u < + \infty $ (by Fatou lemma). Also, $ V_i $ weakly converge to a nonnegative function $ V $ in $ L^{\infty} $. The function $ u $ is in $ W_0^{1, q}(\Omega) $ solution of :

\begin{displaymath}   \left \{ \begin {split} 
      -\Delta u  & = Ve^{u} \in L^1(\Omega)  \,\, &&\text{in} \!\!&&\Omega \subset {\mathbb R}^2, \\
                  u  & = 0 \,\,                                    && \text{in} \!\!&&\partial \Omega.               
\end {split}\right.
\end{displaymath}

According to the corollary 1 of Brezis-Merle result, see [6],   we have $ e^{k u }\in L^1(\Omega), k >1 $. By the elliptic estimates, we have $ u \in W^{2, k}(\Omega)\cap C^{1, \epsilon}(\bar \Omega)$.

We denote by $ f \cdot g $ the inner product of any two vectors $ f $ and $ g $ of $ {\mathbb R}^2 $.

We can write,

\be -\Delta ((u_i-u) \tilde \eta_{\epsilon})= (V_i e^{u_i} -Ve^u)\tilde \eta_{\epsilon} -2\nabla (u_i- u)\cdot \nabla \tilde \eta_{\epsilon} +\epsilon_i (\nabla u_i \cdot x) \tilde \eta_{\epsilon}. \label{(1)}\ee

We use the interior esimate of Brezis-Merle, see [6],

\underbar {\it Step 1:} Estimate of the integral of the first term of the right hand side of $ (\ref{(1)}) $.

We use the Green formula between $ \tilde \eta_{\epsilon} $ and $ u $, we obtain,

\be \int_{\Omega} Ve^u \tilde \eta_{\epsilon} dx =\int_{\partial \Omega} \partial_{\nu} u \eta_{\epsilon} \leq C\epsilon= O(\epsilon) \label{(2)}\ee

We have,

\begin{displaymath}   \left \{ \begin {split} 
      -\Delta u_i- \epsilon_i\nabla u_i\cdot x &= V_ie^{u_i}  \,\, &&\text{in} \!\!&&\Omega \subset {\mathbb R}^2, \\
                                                            u & = 0 \,\,                                     && \text{in} \!\!&&\partial \Omega.               
\end {split}\right.
\end{displaymath}

We use the Green formula between $ u_i $ and $ \tilde \eta_{\epsilon} $ to have:

$$ \int_{\Omega} V_i e^{u_i} \tilde \eta_{\epsilon} dx = \int_{\partial \Omega} \partial_{\nu} u_i \eta_{\epsilon} d\sigma-\epsilon_i\int_{\Omega} (\nabla u_i \cdot x) \tilde \eta_{\epsilon} = $$
\be = \int_{\partial \Omega} \partial_{\nu} u_i \eta_{\epsilon} d\sigma + o(1) \to \mu(\eta_{\epsilon}) \leq \mu(J'_{\epsilon}) \leq 4  \pi - \epsilon_0, \,\,\, \epsilon_0 >0 \label{(3)}\ee

From $ (\ref{(2)}) $ and $ (\ref{(3)}) $ we have for all $ \epsilon >0 $ there is $ i_0 $ such that, for $ i \geq i_0 $,

\be \int_{\Omega} |(V_ie^{u_i}-Ve^u) \tilde \eta_{\epsilon}| dx \leq 4 \pi -\epsilon_0 + C\epsilon \label{(4)}\ee

\underbar { Step 2.1:} Estimate of integral of the second term of the right hand side of $ (\ref{(1)}) $.

Let $ \Sigma_{\epsilon} = \{ x \in \Omega, d(x, \partial \Omega) = \epsilon^3 \} $ and $ \Omega_{\epsilon^3} = \{ x \in \Omega, d(x, \partial \Omega) \geq \epsilon^3 \} $, $ \epsilon > 0 $. Then, for $ \epsilon $ small enough, $ \Sigma_{\epsilon} $ is an hypersurface.

The measure of $ \Omega-\Omega_{\epsilon^3} $ is $ k_2\epsilon^3 \leq meas(\Omega-\Omega_{\epsilon^3}) =\mu_L (\Omega-\Omega_{\epsilon^3}) \leq k_1 \epsilon^3 $.

{\bf Remark}: for the unit ball $ {\bar B(0,1)} $, our new manifold is $ {\bar B(0, 1-\epsilon^3)} $.

(Proof of this fact; let's consider $ d(x, \partial \Omega) = d(x, z_0), z_0 \in \partial \Omega $, this imply that $ (d(x,z_0))^2 \leq (d(x, z))^2 $ for all $ z \in \partial \Omega $ which it is equivalent to $ (z-z_0) \cdot (2x-z-z_0) \leq 0 $ for all $ z \in \partial \Omega $, let's consider a chart around $ z_0 $ and $ \gamma (t) $ a curve in $ \partial \Omega $, we have;

$ (\gamma (t)-\gamma(t_0) \cdot (2x-\gamma(t)-\gamma(t_0)) \leq 0 $ if we divide by $ (t-t_0) $ (with the sign and tend $ t $ to $ t_0 $), we have $ \gamma'(t_0) \cdot (x-\gamma(t_0)) = 0 $, this imply that $ x= z_0-s \nu_0 $ where $ \nu_0 $ is the outward normal of $ \partial \Omega $ at $ z_0 $))

With this fact, we can say that $ S= \{ x, d(x, \partial \Omega) \leq \epsilon \}= \{ x= z_0- s \nu_{z_0}, z_0 \in \partial \Omega, \,\, -\epsilon \leq s \leq \epsilon \} $. It  is sufficient to work on  $ \partial \Omega $. Let's consider a charts $ (z, D=B(z, 4 \epsilon_z), \gamma_z) $ with $ z \in \partial \Omega $ such that $ \cup_z B(z, \epsilon_z) $ is  cover of $ \partial \Omega $ .  One can extract a finite cover $ (B(z_k, \epsilon_k)), k =1, ..., m $, by the area formula the measure of $ S \cap B(z_k, \epsilon_k) $ is less than a $ k\epsilon $ (a $ \epsilon $-rectangle).  For the reverse inequality, it is sufficient to consider one chart around one point of the boundary).

We write,

\be \int_{\Omega} |\nabla ( u_i -u) \cdot \nabla \tilde \eta_{\epsilon} | dx =
\int_{\Omega_{\epsilon^3}} |\nabla (u_i -u) \cdot \nabla \tilde \eta_{\epsilon} | dx + \int_{\Omega - \Omega_{\epsilon^3}} |\nabla (u_i-u) \cdot \nabla \tilde \eta_{\epsilon} | dx. \label{(5)}\ee

\underbar {\it Step 2.1.1:} Estimate of $ \int_{\Omega - \Omega_{\epsilon^3}} |\nabla (u_i-u) \cdot \nabla \tilde \eta_{\epsilon}| dx $.

First, we know from the elliptic estimates that  $ ||\nabla \tilde \eta_{\epsilon}||_{L^{\infty}} \leq C_1 /\epsilon^2 $, $ C_1 $ depends on $ \Omega $

We know that $ (|\nabla u_i|)_i $ is bounded in $ L^q, 1< q < 2 $, we can extract  from this sequence a subsequence which converge weakly to $ h \in L^q $. But, we know that we have locally the uniform convergence to $ |\nabla u| $ (by the Brezis-Merle's theorem), then, $ h= |\nabla u| $ a.e. Let $ q' $ be the conjugate of $ q $.

We have, $  \forall f \in L^{q'}(\Omega)$

$$ \int_{\Omega} |\nabla u_i| f dx \to \int_{\Omega} |\nabla u| f dx $$

If we take $ f= 1_{\Omega-\Omega_{\epsilon^3}} $, we have:

$$ {\rm for } \,\, \epsilon >0 \,\, \exists \,\, i_1 = i_1(\epsilon) \in {\mathbb N}, \,\,\, i \geq  i_1,  \,\, \int_{\Omega-\Omega_{\epsilon^3}} |\nabla u_i| \leq \int_{\Omega-\Omega_{\epsilon^3}} |\nabla u| + \epsilon^3. $$

Then, for $ i \geq i_1(\epsilon) $,

$$ \int_{\Omega-\Omega_{\epsilon^3}} |\nabla u_i| \leq meas(\Omega-\Omega_{\epsilon^3}) ||\nabla u||_{L^{\infty}} + \epsilon^3= \epsilon^3(k_1 ||\nabla u||_{L^{\infty}} + 1) = O(\epsilon^3). $$

Thus, we obtain,

\be \int_{\Omega - \Omega_{\epsilon^3}} |\nabla (u_i-u) \cdot \nabla \tilde \eta_{\epsilon}| dx  \leq \epsilon C_1(2 k_1 ||\nabla u||_{L^{\infty}} + 1)=O(\epsilon)  \label{(6)}\ee

The constant $ C_1 $ does  not depend on $ \epsilon $ but on $ \Omega $.

\underbar {\it Step 2.1.2:} Estimate of $ \int_{\Omega_{\epsilon^3}} |\nabla (u_i-u) \cdot \nabla \tilde \eta_{\epsilon}| dx $.

We know that, $ \Omega_{\epsilon} \subset \subset \Omega $, and ( because of Brezis-Merle's interior estimates) $ u_i \to u $ in $ C^1(\Omega_{\epsilon^3}) $. We have,

$$ ||\nabla (u_i-u)||_{L^{\infty}(\Omega_{\epsilon^3})} \leq \epsilon^3,\, {\rm for } \,\, i \geq i_3. $$

We write,
 
$$ \int_{\Omega_{\epsilon^3}} |\nabla (u_i-u) \cdot \nabla \tilde \eta_{\epsilon}| dx \leq ||\nabla (u_i-u)||_{L^{\infty}(\Omega_{\epsilon^3})} ||\nabla \tilde \eta_{ \epsilon}||_{L^{\infty}} =C_1\epsilon= O(\epsilon) \,\, {\rm for } \,\, i \geq i_3, $$

For $ \epsilon >0 $, we have for $ i \in {\mathbb N} $, $ i \geq i' $,

\be \int_{\Omega} |\nabla (u_i-u) \cdot \nabla \tilde \eta_{\epsilon}| dx \leq \epsilon C_1(2 k_1 ||\nabla u||_{L^{\infty}} + 2) =O(\epsilon) \label{(7)}\ee

From $ (\ref{(4)}) $ and $ (\ref{(7)}) $, we have, for $ \epsilon >0 $, there is $ i'' $ such that, $ i \geq i'' $,

\be \int_{\Omega} |\Delta [(u_i-u)\tilde \eta_{\epsilon}]|dx \leq 4 \pi-\epsilon_0+  \epsilon 2 C_1(2 k_1 ||\nabla u||_{L^{\infty}} + 2 + C) = 4 \pi-\epsilon_0+ O(\epsilon) \label{(8)}\ee

We choose $ \epsilon >0 $ small enough to have a good estimate of  $ (\ref{(1)}) $.

Indeed, we have:

\begin{displaymath}   \left \{ \begin {split} 
      -\Delta [(u_i-u) \tilde \eta_{\epsilon}]   & = g_{i,\epsilon} \,\, &&\text{in} \!\!&&\Omega \subset {\mathbb R}^2, \\
                     (u_i-u) \tilde \eta_{\epsilon}  & = 0 \,\,                                         && \text{in} \!\!&&\partial \Omega.               
\end {split}\right. 
\end{displaymath}

with $ ||g_{i, \epsilon} ||_{L^1(\Omega)} \leq 4 \pi -\epsilon_0/2. $

We can use Theorem 1 of [6] to conclude that there are $ q \geq \tilde q >1 $ such that:

$$ \int_{V_{\epsilon}(x_0)} e^{\tilde q|u_i-u|} dx \leq \int_{\Omega} e^{q|u_i -u| \tilde \eta_{\epsilon}} dx \leq C(\epsilon,\Omega). $$
 
where, $ V_{\epsilon}(x_0) $ is a neighberhooh of $ x_0 $ in $ \bar \Omega $. Here we have used that in a neighborhood of $ x_0 $  by the elliptic estimates, 
$ 1- C \epsilon \leq \tilde \eta_{\epsilon} \leq 1 $.

Thus, for each $ x_0 \in \partial \Omega - \{ \bar x_1,\ldots, \bar x_m \} $ there is $ \epsilon_0>0, q_0 > 1 $ such that:

$$ \int_{B(x_0, \epsilon_{0})} e^{q_0 u_i} dx \leq C, \,\,\, \forall \,\,\, i. $$

By the elliptic estimate see [14] we have:

$$ || u_i||_{C^{1,\theta}[B(x_0, \epsilon)]} \leq c_3 \,\,\, \forall \,\,\, i. $$

We have proved that, there is a finite number of points $ \bar x_1, \ldots, \bar x_m $ such that the squence $ (u_i)_i  $ is locally uniformly bounded in $ C^{1, \theta}, (\theta >0) $ on $ \bar \Omega - \{ \bar x_1, \ldots , \bar x_m \} $.

\underbar {\it Proof of theorem 1.2:} 

The Pohozaev identity gives :

$$  \int_{\partial \Omega} \dfrac{1}{2} (x \cdot \nu) (\partial_{\nu} u_i)^2 d\sigma + \epsilon \int_{\Omega} (x \cdot \nabla u_i)^2 dx+\int_{\partial \Omega} (x \cdot \nu)V_i e^{u_i} d\sigma =  \int_{\Omega} (x \cdot \nabla V_i +2V_i)e^{u_i} dx $$

We use the boundary condition and the fact that $ \Omega $ is starshaped and  the fact that $ \epsilon >0 $ to have that:

\be \int_{\partial \Omega} (\partial_{\nu} u_i)^2 dx  \leq c_0(b, A, C, \Omega). \label{(9)}\ee

Thus we can use the weak convergence in $ L^2(\partial \Omega) $ to have a subsequence $ \partial_{\nu} u_i $, such that: 

$$ \int_{\partial \Omega} \partial_{\nu} u_i \phi dx  \to  \int_{\partial \Omega} \partial_{\nu} u \phi dx, \,\,\, \forall \,\,\, \phi \in L^2(\partial \Omega), $$

Thus, $ \alpha_j = 0 $, $ j=1, \ldots, N $ and $ (u_i) $ is uniformly bounded.

{\bf Remark 1:} Note that if we assume the open set bounded starshaped and $ V_i $ uniformly Lipschitzian and between two positive constants we can bound, by using the inner normal derivative $ \int_{\Omega} e^{u_i} $. 

{\bf Remark 2:} If we assume the open set bounded starshaped and $ \nabla \log V_i $ uniformly bounded, by the previous Pohozaev identity (we consider the inner normal derivative) one can bound $ \int_{\Omega} V_i e^{u_i} $ uniformly.

{\bf Remark 3:} One can consider the problem on the unit ball and an ellipse. These two problems are differents, because: 

1) if we use a linear transformation, $ (y_1,y_2)=(x_1/a, x_2/b) $, the Laplcian is not invariant under this map.

2) If we use a conformal transformation, by a Riemann theorem, the quantity $ x\cdot\nabla u $ is not invariant under this map.

We can not use, after using those transofmations, the Pohozaev identity.

\section{A counterexample}

We start with the notation of the counterexample of Brezis and Merle.

The domain $\Omega $ is the unit ball centered in $ (1,0) $. 

Let's consider $ z_i $ (obtained by the variational method), such that:

$$ -\Delta z_i -\epsilon_i x \cdot \nabla z_i =-L_{\epsilon_i}(z_i)=f_{\epsilon_i}. $$

With Dirichlet condition. By the regularity theorem we have $ z_i \in C^1(\bar \Omega) $.

We have:

$$ ||f_{\epsilon_i}||_1=4\pi A. $$

Thus by the duality theorem of Stampacchia or Brezis-Strauss, we have:

$$||\nabla z_i||_q \leq C_q,\,\, 1\leq q < 2. $$

We solve:

$$ -\Delta w_i =\epsilon_i x \cdot \nabla z_i, $$

With Dirichlet condition.

By the elliptic estimates, $ w_i \in C^1(\bar \Omega) $ and $ w_i \in C^0(\bar \Omega) $ uniformly.

By the maximum principle we have:

$$ z_i-w_i\equiv u_i. $$

Where $ u_i $ is the function of the counterexemple of Brezis Merle.

We write:

$$ -\Delta z_i-\epsilon_i x \cdot \nabla z_i=f_{\epsilon_i}=V_i e^{z_i}. $$

Thus, we have:

$$ \int_{\Omega} e^{z_i} \leq C_1,\,\, {\rm and} \,\, 0\leq V_i \leq C_2, $$

and,

$$ z_i(a_i)\geq u_i(a_i)-C_3 \to +\infty, \,\, a_i\to O. $$

\end{document}